\documentclass{article}
\usepackage[utf8]{inputenc}
\usepackage{amsmath}
\usepackage{mathtools}
\usepackage{amsfonts}
\usepackage{amsthm}
\usepackage{float}
\usepackage[nobysame]{amsrefs}
\usepackage{amssymb,thmtools}
\usepackage{bm}
\usepackage{tikz}
\usepackage{url}
\usetikzlibrary{positioning}
\usetikzlibrary{arrows.meta}
\usepackage{indentfirst}
\usepackage{hyperref}
\usepackage[capitalise]{cleveref}
\usepackage{blkarray}
\usepackage[margin=1in]{geometry}

\AtEndDocument{\bigskip{\footnotesize%
  J.~Boretsky, \textsc{Department of Mathematics, Harvard University, Cambridge, MA, 02138} \par  
  \textit{E-mail address:} \texttt{jboretsky@math.harvard.edu} \par
  \addvspace{\medskipamount}
}}

\newif\ifptitle
\newif\ifpnumber
\newcounter{para}
\newcommand\ptitle[1]{\par\refstepcounter{para}
{\ifpnumber{\noindent\textcolor{lightgray}{\textbf{\thepara}}\indent}\fi}
{\ifptitle{\textbf{[{#1}]}}\fi}}

\makeatletter
\g@addto@macro\@floatboxreset\centering
\makeatother

\usepackage{microtype}
\usepackage{geometry}
\hypersetup{
	pdfstartview={XYZ null null 1.00}, 
	pdfpagemode=UseNone, 
	colorlinks,
	breaklinks, 
	linkcolor=blue,
	urlcolor=blue, 
	anchorcolor=blue,
	citecolor=blue
}

\definecolor{dark-gray}{gray}{0.35}
\newcommand{\Pl}{Pl{\"u}cker }

\newcommand{\TP}{trFl^{\geq 0}}
\newcommand{\trop}{\overline{\textnormal{trop}}}
\newcommand{\Trop}{\textnormal{Trop}\;}

\newtheorem{definition}{Definition}
\newtheorem{theorem}[definition]{Theorem}
\newtheorem{proposition}[definition]{Proposition}
\newtheorem{cor}[definition]{Corollary}
\newtheorem{example}[definition]{Example}
\newtheorem{lemma}[definition]{Lemma}
\newtheorem*{intro}{Theorem}

\numberwithin{definition}{section}
\numberwithin{theorem}{section}
\numberwithin{proposition}{section}
\numberwithin{cor}{section}
\numberwithin{example}{section}
\numberwithin{lemma}{section}

\title{Positive Tropical Flags and the Positive Tropical Dressian}
\author{Jonathan Boretsky\thanks{We acknowledge the support of the Natural Sciences and Engineering Research Council of Canada (NSERC), [funding reference number 557353-2021]. Cette recherche a été financée par le Conseil de recherches en sciences naturelles et en génie du Canada (CRSNG), [numéro de référence 557353-2021].}}

\begin{document}

\maketitle
\begin{abstract}
    We study the totally non-negative part of the complete flag variety and of its tropicalization. We start by showing that Lusztig's notion of non-negative complete flag variety coincides with the flags in the complete flag variety  which have non-negative Pl{\"u}cker coordinates. This mirrors the characterization of the totally non-negative Grassmannian as those points in the Grassmannian with all non-negative Pl{\"u}cker coordinates. We then study the tropical complete flag variety and complete flag Dressian, which are two tropical versions of the complete flag variety, capturing realizable and abstract flags of tropical linear spaces, respectively. The complete flag Dressian properly contains the tropical complete flag variety. However, we show that the totally non-negative parts of these spaces coincide.
\end{abstract}

\section{Introduction}

\ptitle{Motivation for this work} The \textit{Grassmannian} of $k$ planes in $n$ space describes $k$ dimensional linear subspaces in $n$ dimensional space. It is an algebraic variety cut out by the \emph{\Pl relations}. We can \textit{tropicalize} these relations to obtain the \emph{tropical \Pl relations}. The set of points satisfying the tropical \Pl relations, called the \textit{Dressian}, is the parameter space of abstract tropical linear spaces \cite{Spe}. The set of points satisfying the tropicalizations of all polynomials in the ideal generated by the \Pl relations, called the \textit{tropical Grassmannian}, is the parameter space of realizable tropical linear spaces \cite{HKT}. In general, the Dressian properly contains the tropical Grassmannian (see, for instance, \cite{HJM}). However, in \cite{SW2}, it is shown that if we restrict to positive solutions, for an appropriate notion of positivity, the situation is simpler: the \textit{positive Dressian} equals the \textit{positive tropical Grassmannian}. More explicitly, this means that a positive solution to the tropicalizations of the \Pl relations is also a positive solution to the tropicalization of any polynomial in the ideal generated by the \Pl relations. Our goal is to generalize this fact to the setting of the complete flag variety.

\ptitle{Introduction of the major objects in this paper} The \textit{complete flag variety}, $Fl_n$, is the set of complete flags of linear subspaces $\{0\}=V_0\subsetneq V_1\subsetneq\cdots \subsetneq V_n=\mathbb{R}^n$. Any point of this variety is determined by a set of coordinates called its \textit{\Pl coordinates}. These are cut out by the \textit{incidence-\Pl relations}, a set of polynomials which extends the \Pl relations, which generate an ideal called the \textit{incidence-\Pl ideal}. We consider the set of points satisfying the tropicalizations of the incidence-\Pl relations, called the \textit{complete flag Dressian}, $FlDr_n$, and the set of points satisfying the tropicalizations of all polynomials in the incidence-\Pl ideal, called the \textit{tropical complete flag variety}, $TrFl_n$. These parameterize abstract flags of tropical linear spaces and realizable flags of tropical linear spaces, respectively \cite{BEZ}.  

\ptitle{First main theorem} The tropical spaces $FlDr_n$ and $TrFl_n$ are generally different \cite{BEZ}. Motivated by the example of the tropical Grassmannian, we will investigate the totally non-negative (TNN) parts of these spaces. We define the \textit{totally non-negative complete flag Dressian} to be the set of simultaneous positive solutions to the tropicalizations of the incidence-\Pl relations and the \textit{totally non-negative tropical complete flag variety} to be the set of simultaneous positive solutions to the tropicalizations of all the polynomials in the incidence-\Pl ideal. Our main result, \cref{mainthm}, says the following: 

\begin{intro}
The TNN tropical complete flag variety, $TrFl_n^{\geq{0}}$, equals the TNN complete flag Dressian, $FlDr_n^{\geq 0}$.
\end{intro}

\ptitle{Second main theorem} A number of authors, among them \cite{Rie}, \cite{TW}, \cite{Lam} and \cite{Lus3}, have proven that the TNN Grassmannian, in the sense of Lusztig \cite{Lut}, consists precisely of points in the Grassmannian where each \Pl coordinate is non-negative. We extend this result to the setting of the complete flag variety. Specifically, in proving theorem \cref{mainthm}, we will need to carefully study the \textit{totally non-negative complete flag variety}, denoted $Fl_n^{\geq0}$. A construction based on the parameterization of $Fl_n^{\geq 0}$ by Marsh and Rietsch \cite{MR} will allow us to understand explicitly the \Pl coordinates $\{P_I(F)\}_{I\subset [n]}$ of an arbitrary flag $F$ in $Fl_n^{\geq0}$. In \cref{nonnegPluckCoords}, we show:

\begin{intro}
The TNN complete flag variety $Fl_n^{\geq 0}$ equals the set $\{F\in Fl_n|\;P_I(F)\geq 0\; \forall \; I\subset [n]\}$.
\end{intro}

We have learned recently that this result has been independently proven in \cite{BK}, where they show moreover that the only partial flag variety for which this theorem holds are those where the dimensions of the constituent subspaces are consecutive integers. This includes $FL_n^{\geq 0}$, with constituent dimensions $\{1,2,\cdots,n\}$, and the TNN Grassmannian of $k$ planes in $n$ space, with constituent dimension $\{k\}$.

\ptitle{Tangential result of independent interest} We also introduce an alternative description of $FlDr_n$ as the set of points which are common solutions to the three term tropical incidence-\Pl relations and whose supports form flag matroids. This is analogous to the description of the Dressian as the set of common solutions of all the three-term tropical \Pl relations \cite{Mur2}*{Theorem 5.2.25}. While this result is independent of the rest of this abstract, it motivates the introduction of a related TNN tropical space which is used in the proof of \cref{mainthm}.

\ptitle{Structure} The structure of this extended abstract is as follows: In section 2, we introduce the TNN complete flag variety. In section 3, we give a parametrization of this space and study its \Pl coordinates. In section 4, we introduce three tropicalizations of the complete flag variety and demonstrate that two of them are in fact the same. We then focus in on the TNN parts of these three tropical spaces and demonstrate that they are all equal. 

\ptitle{Acknowledgements}\textbf{Acknowledgements:} I would like to thank my supervisor Lauren Williams for introducing me to the tropical flag variety and for many helpful conversations as this work developed. I would also like to thank both Melissa Sherman-Bennett and Mario Sanchez for helpful conversations, examples and references.

 \section{The Totally Non-Negative Complete Flag Variety}
 
\ptitle{Motivation for definition of incidence-\Pl relations}

 \begin{definition}
 The \textbf{complete flag variety} $Fl_n$ is the collection of all \textbf{complete flags} in $\mathbb{R}^n$, which are collections $(V_i)_{i=0}^n$ of linear subspaces satisfying $\{0\}=V_0\subsetneq V_1 \subsetneq \cdots \subsetneq V_n=\mathbb{R}^n$.  
 \end{definition}

 We first observe that $Fl_n$ is a multi-projective variety. We can represent a flag $(V_i)_{i=1}^n$ by a full rank $n$ by $n$ matrix $M$ such that $V_i$ equals the span of the topmost $i$ rows of $M$. Let $GL_n$ be the group of invertible $n$ by $n$ matrices and $SL(n, \mathbb{R})$ be the \textit{special linear group} of real matrices with determinant $1$. Let $B_-$ be the \textit{Borel subgroup} of $GL_n$ consisting of lower triangular matrices. One can check that two matrices $M$ and $M'$ represent the same flag if and only if they are related by left multiplication by some $B\in B_{-}$.  Thus, we can think of the complete flag variety as $Fl_n=\{B_-u|u\in SL(n,\mathbb{R})\}$, where a flag in $Fl_n$ represented by a matrix $u$ is identified with the set $B_-u$.
 
 For $I\subset[n]=\{1,\cdots,n\}$ and $M$ an $n$ by $n$ matrix, the \textit{\Pl coordinate} (or, alternatively, \textit{flag minor}) $P_I(M)$ is the determinant of the submatrix of $M$ in rows $\{1,2,\cdots,|I|\}$ and columns $I$. To any flag $F$, associate the collection of \Pl coordinates $(P_I(F))_{I\subset[n]}$, defined to be the \Pl coordinates of any matrix representative of that flag.  By \cite{MS}*{Proposition 14.2}, this is an embedding of $Fl_n$ in $\mathbb{R}\mathbb{P}^{{n\choose 1}-1}\times\cdots \times \mathbb{R}\mathbb{P}^{{n\choose n-1}-1}$.  The \Pl coordinates of flags in $Fl_n$ are cut out by multi-homogeneous polynomials, as shown in the following definition and theorem.  Note that we will often use shorthand notation such as $(S\setminus ab)\cup cd$ in place of $\left(S\setminus \{a,b\}\right)\cup \{c,d\}$.

 \ptitle{Incidence-\Pl relations}\begin{definition}[\cite{Ful}]
  Consider $\mathbb{R}\mathbb{P}^{{n\choose 1}-1}\times\cdots \times \mathbb{R}\mathbb{P}^{{n\choose n-1}-1}$, with coordinates indexed by proper subsets of $[n]$. For $1\leq r\leq s\leq n$, the \textbf{Incidence-\Pl relations} for indices of size $r$ and $s$ are 
  \begin{equation}
      \mathfrak{P}_{r,s;n}=\left\{\sum_{j\in J\setminus I}sign(j,I,J)P_{I\cup j}P_{J\setminus j }\left| I \in {n\choose r-1},\: J\in {n\choose s+1} \right.\right\},
  \end{equation}
 where $sign(j,I,J)=(-1)^{|\{k\in J|k<j\}|+|\{i\in I|j<i\}|}$.
 
The full set of incidence-\Pl relations is given by $\mathfrak{P}_{IP;n}=\bigcup_{1\leq r\leq s\leq n} \mathfrak{P}_{r,s;n}$. The ideal generated by $\mathfrak{P}_{IP;n}$, denoted $I_{{IP;n}}$, is called the \textbf{incidence-\Pl ideal}.

 \end{definition}
 
\textit{Remark:} Note that the above definition allows for the option of $r=s$. The incidence-\Pl relations for which $r=s$ are called the (Grassmann) \textbf{\Pl relations}. 
  
\ptitle{The complete flag variety from incidence \Pl ideal}

\begin{theorem}[\cite{Ful}*{Section 9, Proposition 1} and discussion following its proof]
 Let $P\in\mathbb{R}\mathbb{P}^{{n\choose 1}-1}\times\cdots \times \mathbb{R}\mathbb{P}^{{n\choose n-1}-1}$. Then $P=P(F)$ for some $F\in Fl_n$ if and only if $P$ satisfies the incidence-\Pl relations $\mathfrak{P}_{IP;n}$. 
 \end{theorem}

In particular, this means the incidence-\Pl relations are precisely the relations between the topmost minors of a generic full rank matrix.

Lusztig introduced the notion of non-negativity for flag varieties. We outline here the definition of the \textit{totally non-negative complete flag variety}, following \cite{Lus2}. We work in type A and so the appropriate simplifications will be made in presenting the definition. Let $s_i$ be the transposition $(i, i+1)$ in the symmetric group $S_n$ and let $w_0$ be the longest permutation in $S_n$. For $1\leq k < n$, let $x_k(a)$ be the $n$ by $n$ matrix which is the identity matrix with an $a$ added in row $k$ of column $k+1$. Explicitly,
\begin{equation*}
    x_k(a)=\begin{blockarray}{ccccccc}
& &  & k & k+1 & & \\
\begin{block}{c(cccccc)}
  & 1 &  &  & & & \\
  & & \ddots &  & &  & \\
  k& &  & 1 & a &  &  \\
  k+1& &  &0  & 1 &  &  \\
  & &  &  &  & \ddots  & \\
  &&&&&&1\\
\end{block}
\end{blockarray} \hspace{5pt},
\end{equation*}
\noindent where unmarked off-diagonal matrix entries are $0$.

\begin{definition}[\cite{Lut}]
 Let $N={n\choose2}$. Pick $\left(i_1,i_2,\cdots, i_N\right)$ such that $s_{i_1}\cdots s_{i_{N}}=w_0$. Then let 
 \begin{equation*}U_{>0}^{+}=\left\{\left.x_{i_1}(a_1)\cdots x_{i_N}\left(a_{N}\right)\right|a_i\in \mathbb{R}_{>0} \;\forall \; i\right\}
 \end{equation*}
\end{definition}

This definition is independent of the choice of sequence $\left(i_1,\cdots, i_{N}\right)$. 

\begin{definition}[\cite{Lut}]\label{lusdef}
 Let $\mathcal{B}_{>0}=\{B_{-}u|u\in U^+_{>0}\}\subset Fl_n$. The \textbf{totally non-negative complete flag variety} (of type A), $Fl_n^{\geq0}$, is the closure of $\mathcal{B}_{>0}$. 
\end{definition}

\section{Parametrization of the TNN  Complete Flag Variety}

\subsection{The Marsh-Rietsch Parametrization}

\ptitle{Introduction to MR parameterization} As shown by Rietsch \cite{Rie}, $Fl_n^{\geq 0}$ is a cell complex, whose cells $\mathcal{R}_{v,w}^{>0}$ are indexed by pairs of permutations $v\leq w$ in the Bruhat order on $S_n$. Each such $\mathcal{R}_{v,w}^{>0}$ is given an explicit parameterization in \cite{MR}. We will describe this parameterization here, making some choices that in principle are arbitrary but will be convenient for our purposes, and invite the reader to look at the above references for full generalities.  

\ptitle{Useful notation}

\ptitle{Positive Distinguished subexpressions, with an intuitive description} Any permutation $w$ in $S_n$ can be written as a product of simple transpositions $s_i$, called an \textit{expression} for $w$. The \textit{length} of $w$, $\ell(w)$, is the fewest number of transpositions in any expression for $w$. An expression for $w$ consisting of $\ell(w)$ transpositions is called \textit{reduced}. Let $\bm{w}=s_{i_1}s_{i_2}\cdots s_{i_k}$ be a reduced expression for $w$. If $v\leq w$ in the Bruhat order, then there is a reduced \textit{subexpression} $\bm{v}=s_{i_{j_1}}s_{i_{j_2}}\cdots s_{i_{j_m}}$ for $v$ in $\bm{w}$, where $1\leq j_1<j_2<\cdots <j_m\leq k$. We will be interested in a special choice of subexpression which is called the \textit{positive distinguished subexpression}. Intuitively, this can be thought of as the leftmost subexpression.

\begin{definition}\label{negative} Let $v\leq w$. Choose a a reduced expression $\bm{w}=s_{i_1}s_{i_2}\cdots s_{i_k}$ for $w$ and a subexpression $\bm{v}=s_{i_{j_1}}\cdots s_{i_{j_m}}$ for $v$ in $\bm{w}$. Then $\bm{v}$ is a \textbf{positive distinguished subexpression} if:

\begin{enumerate}
    \item It is a reduced expression for $v$.
    \item Whenever $\ell(s_{i_p}s_{i_{j_r}}\cdots s_{i_{j_m}})<\ell(s_{i_{j_r}}\cdots s_{i_{j_m}})$ for $j_{r-1}\leq p< j_{r}$, we have $p=j_{r-1}$. 
\end{enumerate}  
\end{definition}

\begin{lemma}[\cite{MR}*{Lemma 3.5}]
For every $v\leq w$, and every reduced expression $\bm{w}$ of $w$, there is a unique positive distinguished subexpression for $v$ in $\bm{w}$.
\end{lemma}

\ptitle{Example of positive distinguished subexpression}
\begin{example}
Let $n=4$. Set $\bm{w}=s_1s_2s_3s_1s_2s_1$ and $v=s_1s_2s_1$. The leftmost subexpression for $v$ in $\bm{w}$ is $j_1=1$, $j_2=2$ and $j_3=4$. Indeed, one can verify that this choice satisfies the definition.
\end{example}

For $1\leq k<n$, let $\dot{s}_k$ be the $n$ by $n$ identity matrix with the $2\times 2$ submatrix in rows $\{k,k+1\}$ and columns $\{k,k+1\}$ replaced by the matrix $\begin{pmatrix}
0&1\\
-1&0
\end{pmatrix}$. Explicitly, 
\begin{equation*}
    \dot{s}_k=\begin{blockarray}{ccccccc}
& &  & k & k+1 & & \\
\begin{block}{c(cccccc)}
  & 1 &  &  & & & \\
  & & \ddots &  & &  & \\
  k& &  & 0 & 1 &  &  \\
  k+1& &  &-1  & 0 &  &  \\
  & &  &  &  & \ddots  & \\
  &&&&&&1\\
\end{block}
\end{blockarray}\hspace{5pt},
\end{equation*}
\noindent where unmarked off-diagonal matrix entries are $0$.

We will describe each cell of $Fl_n^{\geq 0}$ as a product of matrices of the form $x_k$ and $\dot{s}_{k'}$.

\begin{definition}
Fix $v\leq w$ in the Bruhat order. Fix a vector $\bm{a}\in \mathbb{R}^{\ell(w)-\ell(v)}$. Consider the reduced expression $\bm{w_0}=(s_1s_2\cdots s_{n-1})(s_1s_2\cdots s_{n-2})(\cdots)(s_1s_2)(s_1)$ for $w_0$, the longest permutation in the Bruhat order in $S_n$\footnote{This choice of expression is arbitrary in the context of the Marsh-Rietsch parameterization, but plays an important role in the proofs underlying later results in this abstract.}. Choose the positive distinguished subexpression $\bm{w}$ for $w$ in $\bm{w_0}$, and the positive distinguished subexpression $\bm{v}$ for $v$ in $\bm{w}$, and write them as $\bm{w}= s_{i_1}\cdots s_{i_k}$ and $\bm{v}=s_{i_{j_1}}\cdots s_{i_{j_m}}$, respectively. Let $J=\{j\;|\;j=j_t \textnormal{ for some }t\}$. In other words, $J$ are those indices which correspond to transpositions that are used in $\bm{v}$. Then set
\begin{equation*}
    M_{v,w}(\bm{a})\coloneqq M_1\cdots M_k,\hspace{3mm} \textnormal{  where  }\hspace{3mm} M_j=\begin{cases}\dot{s}_{i_j}, & j\in J\\ x_{i_j}(a_{j}),&j\notin J\end{cases}\hspace{1mm} .
\end{equation*}
\end{definition}

 \begin{theorem}[Marsh-Rietsch Parametrization, \cite{MR}] \label{marshrietsch}

Each cell $\mathcal{R}_{v,w}^{>0}$ of $Fl_{n}^{\geq 0}$ can be parameterized as 

\begin{equation*}
    \mathcal{R}_{v,w}^{>0}=\left\{M_{v,w}(\bm{a})\left|\bm{a}\in\mathbb{R}_{>0}^{\ell(w)-\ell(v)}\right.\right\}
\end{equation*}

In particular, each flag $F\in Fl_n^{\geq 0}$ is uniquely represented in some unique $\mathcal{R}_{v,w}^{> 0}$. Moreover, each $\mathcal{R}_{v,w}^{> 0}$ is a cell, meaning it is homeomorphic to an open ball. \end{theorem}

\begin{example}\label{toymodel}
Let $n=4$, $\bm{w}=s_1s_3s_2s_1$ and $v=s_2$. The positive distinguished subexpression for $v$ in $w$ is the subexpression where $j_1=3$, so $J=\{3\}$. Thus, $M_1=x_1(a_1)$, $M_2=x_3(a_2)$, $M_3=\dot{s}_2$ and $M_4=x_1(a_3)$. The cell of the non-negative flag variety corresponding to $v\leq w$ is represented by matrices of the form 

\begin{equation*}
    M=M_1M_2M_3M_4 
    =\begin{pmatrix}
    1&a_3&a_1&0\\
    0&0&1&0\\
    0&-1&0&a_2 \\
    0&0&0&1
    \end{pmatrix},
\end{equation*}
where the $a_i$ range over all positive real numbers. Note that the \Pl coordinates of $M$ are all non-negative, as \cref{nonnegPluckCoords} will show must be true for any cell of $Fl_n^{\geq 0}$. 
 \end{example}
 
 \ptitle{Useful results for later}
 We now give a useful property of the cells $\mathcal{R}_{v,w}^{> 0 }$.
 
 \begin{lemma}\label{cellcoords}
 Each cell $\mathcal{R}_{v,w}^{> 0}$ of $Fl_n^{\geq 0}$ consists entirely of flags for which some fixed collection of \Pl coordinates is strictly positive and the rest are $0$.   
 \end{lemma}

\subsection{Extremal Non-Zero \Pl Coordinates}

 \ptitle{Outline/goal of the section} We define a special subset of the \Pl coordinates of a flag which we call \textit{extremal non-zero \Pl coordinates}. The set of indices of the extremal non-zero \Pl coordinates of a flag in $Fl_n^{\geq 0}$ will depend only on which cell $\mathcal{R}_{v,w}^{>0}$ that flag lies in. Further, we will show that in any given cell of $Fl_n^{\geq 0}$, the extremal non-zero \Pl coordinates determine all of the other \Pl coordinates.

\ptitle{Description and intuition of $\Xi$} \ptitle{Description and intuition of $\Xi$} For any $1\leq k<n$ and any $P\in \mathbb{R}\mathbb{P}^{{n\choose 1}-1}\times\cdots\times \mathbb{R}\mathbb{P}^{{n\choose n-1}-1}$, we define a map $\Xi_P: {[n]\choose k}\rightarrow {[n]\choose k}$.  Intuitively, when applied to the index of a non-zero \Pl coordinate $I$, this map finds the largest member of $I$ that can be increased without making the corresponding \Pl coordinate $0$ and increases it maximally. Explicitly, for $I$ such that $P_I\neq 0$, define $b=\max_{i\in I}\left\{i\;|\;\exists\; j,  \;i<j\notin I, \; P_{(I\setminus i)\cup j}\neq 0\right\}$, if that set is non-empty. Otherwise, say $b$ does not exist. If $b$ exists, define $a=\max_{j\notin I}\left\{j\;|\; P_{(I\setminus b)\cup j}\neq 0\right\}$. Then, 
\begin{align*}
    \Xi_P(I)=
    \begin{cases}
    (I\setminus b)\cup a &\;\;\textnormal{if } I \textnormal{ is the index of a non-zero \Pl coordinate and $b$ exists,}\\
    I &\;\;\textnormal{otherwise.}  
    \end{cases}
\end{align*}

Note that the indices of non-zero \Pl coordinates with index of some fixed size can be seen as the bases of a matroid. In this light, $\Xi_P$ acts by basis exchange. Also note that, by \cref{cellcoords}, for a TNN flag $F$, the map $\Xi_{P(F)}$ depends only on the cell $\mathcal{R}_{v,w}$ in which $F$ lies. 

The extremal non-zero \Pl coordinates will be indexed by certain $\Xi$ orbits. To properly define them, we first need a preliminary result on matroids:

 \begin{definition}
  The \textbf{Gale order} on subsets of $[n]$ of size $k$ is a partial order such that, if $I=\{i_1<\cdots< i_k\}$ and $J=\{j_1<\cdots< j_k\}$, then we say $I\leq J$ if $i_r\leq j_r$ for every $r\in[k]$.
 \end{definition}

 \begin{lemma}[\cite{Coxeter}*{Theorem 1.3.1}]  
 Any matroid has a unique Gale minimal basis and a unique Gale maximal basis. 
 \end{lemma}
 
Note that the Gale minimal and maximal bases referenced in the previous lemma must simply be the lexicographically minimal and maximal bases, respectively.

 \begin{definition}\label{extremal}
  Given a set of \Pl coordinates $\{P_I\}$ of a flag, let $I_k$ be the Gale minimal index of size $k$ such that $P_{I_k}\neq 0$.
  The set of indices of the \textbf{extremal non-zero \Pl coordinates} (referred to as \textbf{extremal indices}) 
  of a point $P$ in $\mathbb{R}\mathbb{P}^{{n\choose 1}-1}\times\cdots \times \mathbb{R}\mathbb{P}^{{n\choose n-1}-1}$ is the set consisting of those indices which are in the $\Xi_P$ orbit of $I_k$ for some $k\in [n-1]$.
 
  \end{definition}
  
 If $F$ is a TNN flag, the extremal indices of the \Pl coordinates $P(F)$ depend only on the cell $\mathcal{R}_{v,w}^{>0}$ in which $F$ lies, since $\Xi_{P(F)}$ depends only on the cell in which $F$ lies.
 
 If we have a collection of \Pl coordinates $\{P_I\}$ such that $P_I\geq 0$ for all $I\subset {[n]\choose k}$, the indices of the non-zero \Pl coordinates of $P$ of fixed size form not just the bases of a matroid, but the \textit{bases of a positroid.} For background on positroids, see \cite{Oh}. Many of the useful properties of the extremal \Pl coordinates only hold in this context of positroids.

 \ptitle{Example of extremal indices and \Pl coordinates}
   \begin{example}
 Let $a,b,c,d,e,f,g\in\mathbb{R}_{>0}$ and consider  
 \begin{equation*}
     M=\begin{pmatrix}
    1&a+e+g&ab+af+ef&abc&abcd\\
    0&1&b+f&bc&bcd\\
    0&0&1&c&cd \\
    0&0&0&1&d\\
    0&0&0&0&1
    \end{pmatrix}.
 \end{equation*}
 
 The minors of this matrix are non-negative and specifically, focusing on indices of size $2$, we can see that they are all positive except for $P_{45}=0$. Thus, the non-zero \Pl coordinate with Gale minimal index of size $2$ is $P_{12}$. Then, $\Xi_{P(M)}(12)=15$, replacing the $2$ with a $5$. Next, $\Xi_{P(M)}(15)=35$, replacing the $1$ with a $3$. Thus, $P_{12}=1$, $P_{15}=bcd$ and $P_{35}=bcdef$ are the extremal non-zero \Pl coordinates of size $2$ of this flag. 
 \end{example}

The next theorem highlights the importance of the extremal non-zero \Pl coordinates.
 
 \begin{theorem}\label{3termgenerate}
 For any flag $F$ with non-negative \Pl coordinates, the extremal non-zero \Pl coordinates of $F$ uniquely determine the other non-zero \Pl coordinates of $F$ by three-term incidence-\Pl relations. 
 \end{theorem}

 \subsection{\Pl Coordinates of the TNN Flag Variety}
 
 \ptitle{Motivation for this subsection} Now, given a set of extremal non-zero \Pl coordinates for a flag lying in $\mathcal{R}_{v,w}^{> 0}$, we want to understand how to construct a set of parameters $a_i$ for which \cref{marshrietsch} yields a matrix agreeing with those coordinates. 

\ptitle{Establishment of notation and the first main result of this section} 


\begin{theorem}\label{pointtograph}
For any $v\leq w$ with $r=\ell(w)-\ell(v)$, let $\Psi_{v,w}: \mathcal{R}_{v,w}^{> 0}\rightarrow \mathbb{R}^r$ be the map $M_{v,w}(\bm{a})\mapsto \bm{a}$, in the notation of \cref{marshrietsch}. The map $\Psi_{v,w}$ consists of Laurent monomials in the extremal \Pl coordinates. 
\end{theorem}

We can use this theorem to prove the following, which is one of our main results: 

\begin{theorem}\label{nonnegPluckCoords}
The TNN flag variety defined in \cref{lusdef} is precisely the set of flags with non-negative \Pl coordinates. In other words, $Fl_n^{\geq 0}=\{F\in Fl_n|\;P_I(F)\geq 0\; \forall \; I\subset [n]\}$.
\end{theorem}

 It is shown in \cite{KW}*{Lemma 3.10} that any flag in $Fl_n^{\geq 0}$ has non-negative \Pl coordinates. We now outline the strategy used to obtain the converse. 

\begin{definition} 
  A (complete) \textbf{flag matroid} on a ground set $E$ of size $n$ is a sequence of matroids $\bm{\mathcal{M}}=(\mathcal{M}_1,\mathcal{M}_2,\cdots, \mathcal{M}_{n-1})$ on the ground set $E$ with the rank of $\mathcal{M}_i$ equal to $i$, called \textbf{constituent matroids}, such that for any $j<k$,
  
  \begin{itemize}
      \item each basis of $\mathcal{M}_j$ is contained in some basis of $\mathcal{M}_k$.
      \item each basis of $\mathcal{M}_k$ contains some basis of $\mathcal{M}_j$. 
  \end{itemize}
   \end{definition}
  We identify a flag matroid with the collection of bases of its constituent matroids, collectively referred to as the \textit{bases of the flag matroid}. For more details about flag matroids, as well as cryptomorphic definitions, see \cite{Coxeter}. Note that the indices of non-zero \Pl coordinates of an invertible square matrix are easily seen to form a flag matroid. 
  
  \begin{definition}A flag matroid on $[n]$ is \textbf{realizable} if its bases are the non-zero \Pl coordinates of some $F\in Fl_n$. 
  \end{definition}

   We now define two types of flag positroid. The difference between their definitions mirrors the apparent difference between a flag lying in $Fl_n^{\geq 0}$ according to \cref{lusdef} and a flag with non-negative \Pl coordinates.
  
  \begin{definition}
  A \textbf{realizable flag positroid} on $[n]$ is the set of indices of non-zero \Pl coordinates of a flag $F\in Fl_n^{\geq 0}$ (as defined in \cref{lusdef}). A \textbf{synthetic flag positroid} on $[n]$ is the set of indices of non-zero \Pl coordinates of a flag $F$ satisfying $P_I(F)\geq 0$ for all $I\subset [n]$. 
  \end{definition}
  A priori, one may expect that there could be more synthetic flag positroids than realizable flag positroids, but this is not the case.

  \begin{theorem}\label{synthetic}
 The set of synthetic flag positroids on $[n]$ equals the set of realizable flag positroids on $[n]$.
  \end{theorem}
  
  Note that by \cref{cellcoords}, the realizable flag positroid arising from the non-zero \Pl coordinates of a TNN flag only depends on which cell $\mathcal{R}_{v,w}^{>0}$ that flag lies in. Thus, we can associate a cell $\mathcal{R}_{v,w}^{>0}$ to any realizable flag positroid. Let $F$ be a flag whose \Pl coordinates $P$ are all non-negative. Let $\bm{\mathcal{M}}$ be the synthetic (equivalently, realizable) flag positroid which has $I\subset [n]$ as a basis if and only if $P_I> 0$. As above, let $\mathcal{R}_{v,w}^{>0}$ be the cell associated to $\bm{\mathcal{M}}$. To prove \cref{nonnegPluckCoords}, we are left to show that $F\in\mathcal{R}_{v,w}^{>0}$. Extending the domain of the map $\Psi_{v,w}$ defined in \cref{pointtograph}, one can apply it to $F$ and show that $M_{v,w}\left(\Psi_{v,w}(F)\right)$ is a flag in $\mathcal{R}_{v,w}^{>0}$ which has the same extremal \Pl coordinates as $F$. Then, using \cref{3termgenerate}, one may conclude that $F$ itself lies in $\mathcal{R}_{v,w}^{>0}$, completing the proof of \cref{nonnegPluckCoords}.

\section{Tropicalizing the Complete Flag Variety}

\ptitle{Section outline} We now discuss \textit{tropical varieties} and introduce the precise definitions of the \textit{TNN tropical complete flag variety} and the \textit{TNN complete flag Dressian}.

\ptitle{Definition of tropical polynomials and tropical solutions to polynomials} 

\begin{definition}\label{tropsol}
Let $\bm x=(x_1,\cdots,x_n)\in \mathbb{R}^n$ and $\bm{b}=(b_1,\cdots,b_n)\in\mathbb{N}^n$. We will use the notation $\bm {x}^{\bm {b}}=x_1^{{b}_1}\cdots x_n^{{b}_n}$. Let $p=\sum_i \pm a_i\bm{x}^{\bm b_i}$ be a polynomial, where each $a_i>0$ and each $\bm b_i\in \mathbb{N}^n$. We define the \textbf{tropicalization of $p$} by $trop(p)=\min_i \left\{a_i+\bm{x}\cdot \bm {b_i}\right\}$.  We say that a point $\bm{y}\in\mathbb{T}^n\coloneqq (\mathbb{R}\cup\infty)^n$ is a \textbf{solution of the tropicalization of $p$} if

\begin{equation*}
    \min_i \left\{a_i+\bm{y}\cdot \bm {b_i}\right\}= \min_i \left\{a_i+y_1(b_i)_1+\cdots+ y_n(b_i)_n\right\}
\end{equation*}
is achieved at least twice. We further say that a point in $\mathbb{T}^n$ is a \textbf{positive solution of the tropicalization of $p$} if additionally, at least one of the minima comes from a term of $p$ with a $+$ sign, and at least one of the minima comes from a term with a $-$ sign. Equivalently, if we rewrite $p=0$ in the form $\sum_j c_j\bm{x}^{\bm b_j}=\sum_i a_i\bm{x}^{\bm b_i}$ with all $c_j$ and $a_i$ positive, then we want at least one minimum to occur in a term coming from each side of the equality.
\end{definition}

\ptitle{Projective Tropical spaces} The tropical objects we are interested in will live in \textit{projective tropical spaces}, which are spaces that interact nicely with homogeneous polynomials.

\begin{definition}
 \textbf{Projective tropical space}, denoted $\mathbb{T}\mathbb{P}^n$ is given by $\left(\mathbb{T}^{n+1}\setminus (\infty,\cdots \infty)\right)/\sim$ where the equivalence relation is $\bm x \sim \bm y$ if there exists $c\in \mathbb{R}$ such that $x_i=y_i+c$ for all $i\in [n]$. 
\end{definition}

The following is immediate from the definition: 

\begin{proposition}
If $p$ is a homogeneous polynomial, then $\bm x$ is a (positive) solution of $trop(p)$ if and only if $\bm y$ is a (positive) solution of $trop (p)$ for all $\bm y \sim \bm x$. 
\end{proposition}

\ptitle{Example of tropical solutions}

 \ptitle{Definition of tropical prevariety}

\begin{definition}\label{tropideal}
 Given a set of multi-homogeneous polynomials $\mathcal{P}$, each of which is homogeneous with respect to sets of variables of sizes $\{n_i\}_{i=1}^t$, and the ideal $I$ which they generate, we define the following sets in $\mathbb{T}\mathbb{P}^{n_1-1}\times\cdots\times \mathbb{T}\mathbb{P}^{n_t-1}$:
 
 \begin{itemize}
     \item The \textbf{tropical prevariety} $\trop(\mathcal{P})$ or $\trop(I)$ is the set of simultaneous solutions to the tropicalizations of all the polynomials in $\mathcal{P}$ or in $I$, respectively.
     
     \item The \textbf{non-negative tropical prevariety}, $\trop^{\geq 0}(\mathcal P)$ or $\trop^{\geq 0}(I)$, is the set of simultaneous positive solutions of the tropicalizations of all the polynomials in $\mathcal{P}$ or in $I$, respectively.
 \end{itemize}
 
\end{definition}

\ptitle{Intuition behind positive tropical solutions} Solutions of tropicalizations of polynomials can alternatively be described in a way that more clearly explains the term ``positive solution". Let $\mathcal{C}=\bigcup_{n=1}^\infty\mathbb{C}((t^{1/n}))$ be the field of \textit{Puisseux series} over $\mathbb{C}$. A Puisseux series $p(t)\in \mathcal{C}$ has a term with a lowest exponent, say $at^{u}$ with $a\in \mathbb{C}^*$ and $u\in\mathbb{Q}$. In this case, we define $\textnormal{val}(p(t))=u$. Also, we will define the semifield $\mathcal{R}^+$ to be the set of $p(t)$ in $\mathcal{C}$ where the coefficient of $t^{\textnormal{val}(p(t))}$ is in $\mathbb{R}^+$. In fact, $\mathcal{R}^+$ and $\mathcal{C}$ can be thought of as analogous to $\mathbb{R}^+$ and $\mathbb{C}$, respectively. Given an ideal $I\trianglelefteq \mathcal{\mathbb{C}}[x_1,\cdots x_n]$, let $V(I)\subseteq \mathcal{C}^n$ be the variety where all polynomials in $I$ vanish. We define the \textit{positive part} of this variety to be $V^+(I)=V(I)\cap \mathcal{({R}^+)}^n$. 

\begin{proposition}[\cite{SS} Theorem 2.1 and \cite{SW} Proposition 2.2]\label{otherdef}

Let $I$ be an ideal of $\mathbb{C}[x_1,\cdots, x_n]$. Then $\trop(I)=\overline{\textnormal{val}(V(I))}$ and $\trop^{\geq0}(I)=\overline{\textnormal{val}(V^+(I))}$, where $\overline{\textnormal{val}(V(I))}$ and $\overline{\textnormal{val}(V^+(I))}$ are the closures of $\textnormal{val}(V(I))$ and $\textnormal{val}(V^+(I))$, respectively.

\end{proposition}

\ptitle{The tropical flag variety and the Flag Dressian} Having introduced $Fl_n$, we now define two tropical analogues of this space along with their totally non-negative parts. Recall that $\mathfrak{P}_{IP;n}$ is the set of incidence-\Pl relations and $I_{IP;n}$ is the ideal generated by those relations.

\begin{definition}
 We define the \textbf{tropical complete flag variety} to be $trFl_n=\trop(I_{IP;n})$ and the \textbf{totally non-negative tropical complete flag variety} to be
 $\TP_n=\trop^{\geq 0}(I_{IP;n})$. We define the \textbf{complete flag Dressian} to be $FlDr_n=\trop(\mathfrak{P}_{IP;n})$ and the \textbf{totally non-negative complete flag Dressian} to be $FlDr_n^{\geq 0}=\trop^{\geq0}(\mathfrak{P}_{IP;n})$.   
 
\end{definition}

\ptitle{Prelude and motivation of equality of non-negative parts of Dressian and TFV} \cref{mainthm} will show that $\TP_n$ and $FlDr_n^{\geq 0}$ coincide. Note that this is not obvious, since a point in $trFl_n^{\geq 0}$ satisfies more relations than a point in $FlDr_n^{\geq 0}$. In fact, in general, the tropical prevariety of a collection of polynomials will properly contain the tropical prevariety of the ideal those polynomials generate. In the specific case of the complete flag variety, it is shown in \cite{BEZ}*{Example 5.2.4} that for $n\geq 6$, $FlDr_n$ properly contains $trFl_n$. Before getting to \cref{mainthm}, we need to discuss one other type of Dressian.

\begin{definition}
   The \textbf{support} of $P\in\mathbb{T}\mathbb{P}^{{n\choose 1}-1}\times\cdots\times\mathbb{T}\mathbb{P}^{{n\choose n-1}-1}$ is the set $\left\{I\subset [n]\left|P_I\neq \infty\right.\right\}$. 
  \end{definition}
  
  \ptitle{Thee term flag Dressian and why we want it} 
  \begin{definition}
   Let $\mathfrak{P}^{3\textnormal{ term}}_{IP;n}$ be the set of incidence-\Pl relations with precisely three terms. We denote by $FlDr^{3M}_n$ the subset of $\trop(\mathfrak{P}^{3\textnormal{ term}}_{IP;n})$ consisting of points whose supports form a flag matroid. We call this the \textbf{three-term complete flag Dressian}. Similarly, we denote by ${(FlDr^{3M}_n)}^{\geq 0}$ the subset of $\trop^{\geq 0}(\mathfrak{P}^{3\textnormal{ term}}_{IP;n})$ consisting of points whose supports form a flag matroid, which we call the \textbf{totally non-negative three-term complete flag Dressian}.
  \end{definition}
  
 We can write down the following complete flag version of the fact that the Dressian is cut out by tropical three term \Pl relations \cite{Mur2}*{Theorem 5.2.25}. 
 \ptitle{Equality of the flag Dressian and three-term flag Dressian} 
 \begin{theorem}\label{threegen}
The sets $FlDr^{3M}_n$ and $FlDr_n$ are equal.

 \end{theorem}

We now shift our attention to the non-negative parts of the tropical varieties we have introduced. For $v\leq w$ in the Bruhat order with $r=\ell(w)-\ell(v)$, let $\Phi_{v,w}:\mathbb{R}_{>0}^r\rightarrow \mathbb{R}\mathbb{P}^{{n\choose 1}-1}\times\cdots\times\mathbb{R}\mathbb{P}^{{n\choose n-1}-1}$ be the map which takes a collection of $\bm{a}\in \mathbb{R}_{>0}^r$ to the \Pl coordinates of the matrix $M_{v,w}(\bm{a})$, in the notation of \cref{marshrietsch}. Note that by construction, this map consists of a collection of polynomials in the $a_i$, and so we can tropicalize this map, obtaining a map $\Trop \Phi:\mathbb{R}^r\rightarrow \mathbb{T}\mathbb{P}^{{n\choose 1}-1}\times\cdots\times\mathbb{T}\mathbb{P}^{{n\choose n-1}-1}$. The following corollary can be deduced from \cref{nonnegPluckCoords} with little extra work, and then put to use to help prove our second main result.

\begin{cor}\label{tropvw}
Every point in ${(FlDr^{3M}_n)}^{\geq 0}$ lies in the image of $\Trop \Phi_{v,w}$ for some $v\leq w$.
\end{cor}

\ptitle{The main theorem}
\begin{theorem} \label{mainthm}
The following sets are equal:
\begin{enumerate}
    \item The TNN topical flag variety $trFl^{\geq 0}_n$,
    \item The TNN three-term complete flag Dressian ${(FlDr^{3M}_n)}^{\geq 0}$,
    \item The TNN complete flag Dressian $FlDr_n^{\geq 0}$.
\end{enumerate}
\end{theorem}

\bibliographystyle{plain}
\bibliography{bibliography}

\end{document}